# When is the combined load identifiable from a stress-intensity profile? A coupled forward–inverse study on SIFBench finite-element data

Giansalvo Cirrincione[a], Filippo Grassia[a,*]

[a]*Laboratory of Innovative Technologies (LTI, EA 3899), Université de Picardie Jules Verne, Amiens, 80025, France*

ABSTRACT

This work studies the inverse problem of recovering the relative magnitudes of the tension, bending, and bearing loads acting on a crack from its stress-intensity-factor profile along the crack front, using the public SIFBench finite-element data. The central claim is not forensic load recovery on field cases, but a rigorous characterization of *when* the combined load is identifiable at all, together with an estimator that returns calibrated uncertainty precisely in the regimes where it is not. For a known geometry the forward map from loads to profile is exactly linear, and identifiability reduces to a single geometric question: whether the three elementary load profiles are linearly independent as functions along the front. When they are nearly dependent, many different load combinations produce almost the same profile and the inverse problem is ill-posed; the analysis shows that the degree of ill-posedness is controlled by an intrinsic stability margin, not by the conditioning number alone. A single crack-front operator serves both as a structured forward surrogate and as the differentiable map required by a simplex-constrained, set-valued inverse estimator. On the SIFBench corner-crack scenario the empirical behaviour matches the theory: the typical geometry is well posed while a sizable minority is genuinely ill-posed, so a point estimate is reliable on the majority and provably uninformative on the rest. Validation is on controlled synthetic noise; no real fracture cases are used or claimed.

## 1. Introduction

The stress-intensity factor (SIF) governs fatigue crack growth and is the central quantity in damage-tolerant design. For realistic three-dimensional cracks it varies along the crack front, and computing it traditionally requires finite-element analysis that is too costly for the large parameter sweeps needed in design and life prediction. This has motivated data-driven surrogates, and recently the SIFBench benchmark (Gautam, Kirby, Hochhalter and Zhe, 2025), a public, large-scale dataset of front-resolved SIF profiles across thousands of geometries and three reference loads (tension, bending, bearing). SIFBench also establishes baselines — random forests, support-vector regression, feed-forward networks, and a Fourier neural operator — but exclusively for the *forward* task: predict the SIF profile $K(\varphi)$ along the crack front (parameterized by the front angle $\varphi$) from geometry and a single known load.

The complementary and largely unaddressed question is the *inverse* one. Given an observed front profile, which combination of tension, bending and bearing produced it? This is the quantity of interest in failure analysis and structural diagnostics, yet it has not been studied systematically on a public, reproducible dataset, and it is fundamentally harder than the forward map: different load mixes can produce indistinguishable profiles, so the problem is ill-posed in the Hadamard sense for entire families of geometries. An inverse method that returns a single load estimate without saying when that estimate is meaningful is therefore not merely inaccurate but misleading.

The framing of scope is made explicit and deliberately conservative. Genuine forensic load identification on real failed components would require documented field cases with known loading, which are not available here; the present work therefore does not claim forensic results. "Forensic" and "structural diagnostics" appear in this paper only as motivation. The contribution is methodological: a rigorous account of *when* the combined load is recoverable, an estimator that is exact where it is and honestly set-valued where it is not, and a forward surrogate accurate enough to make that analysis apply to learned rather than tabulated profiles. Validation is on controlled synthetic noise, and the theory makes that validation quantitative rather than ad hoc. This also explains why this is one paper and not two. A stand-alone forward surrogate would be incremental — SIFBench has already benchmarked that task. The value is in the coupling: the same crack-front operator that predicts the profile is the differentiable forward map inside the inverse estimator, and the identifiability theory is a statement about exactly the object the operator produces. Forward and inverse close on a single mathematical object, and separating them would lose precisely what is new.

*Corresponding author
filippo.grassia@u-picardie.fr (F. Grassia)
ORCID(s):

*Notation*
Throughout, $K_T$, $K_B$, $K_P$ denote the front profiles produced by unit tension, bending, and bearing loads; a combined load is the convex mixture with coefficients $(\alpha, \beta, \gamma)$ on the simplex; $M(g)$ is the matrix whose columns are the three profiles for a geometry $g$ (defined in §4); and $\kappa_s$ is the simplex-restricted conditioning number of the inverse problem, made precise in Definition 1. These are used informally in the contributions below and formally from §3 onward.

*Contributions*

**C1** A crack-front operator that is competitive with the strongest SIFBench tabular baseline on the twin corner-crack scenario, improves the $K_T$ and $K_P$ channels, enforces front symmetry and positivity by construction, and supplies the full differentiable profile matrix $M(g)$ required by the inverse estimator. §5.

**C2** A reusable simplex-sampled combined-load dataset ($K_{\text{mix}}$) derived from SIFBench's separated $K_T, K_B, K_P$. §3.

**C3** A simplex-constrained inverse estimator and a calibration protocol whose uncertainty is designed to track the identifiability margin: exact in the known-$M$ regime, with full credible-region coverage stated as a falsifiable diagnostic rather than claimed on the learned map. §6.

**C4** An analytical identifiability theory: exact characterization of when $(\alpha, \beta, \gamma)$ is recoverable from $K(\varphi)$. §4.

## 2. Related work

*SIF surrogates and SIFBench*
Classical SIF determination rests on closed-form and weight-function solutions — the Newman–Raju equations for surface and corner cracks (Newman and Raju, 1981, 1983) and universal weight functions (Glinka and Shen, 1991), building on the foundational linear-elastic crack-tip analysis of Irwin (Irwin, 1957). Data-driven SIF prediction subsequently used Gaussian processes, tree ensembles and neural networks on tabulated finite-element results; SIFBench (Gautam et al., 2025) consolidates this line of work into a public benchmark and shows that, on its front-resolved data, simple tabular regressors are competitive with and often better than a Fourier neural operator. This observation is taken as a design signal rather than a defect: it indicates that the per-front functional structure is being discarded by point-wise models and mishandled by a grid-based operator, which is precisely the gap the crack-front operator of §5 targets.

*Neural operators and coordinate encodings*
Operator-learning methods — Fourier neural operators (Li, Kovachki, Azizzadenesheli, Liu, Bhattacharya, Stuart and Anandkumar, 2021), DeepONet (Lu, Jin, Pang, Zhang and Karniadakis, 2021), and the general neural-operator framework (Kovachki, Li, Liu, Azizzadenesheli, Bhattacharya, Stuart and Anandkumar, 2023) — predict functions rather than point values and are natural when the output is a field over a continuous coordinate. Their effectiveness depends on how the query coordinate is encoded; as Fourier-feature analysis shows (Tancik, Srinivasan, Mildenhall, Fridovich-Keil, Raghavan, Singhal, Ramamoorthi, Barron and Ng, 2020), the spectral content of the encoding controls which output frequencies are representable. The architecture used here is a transformer (Vaswani, Shazeer, Parmar, Uszkoreit, Jones, Gomez, Kaiser and Polosukhin, 2017) over the front coordinate with a Fourier feature encoding (learnable frequencies, though on smooth SIF profiles a fixed spectrum suffices; §5), specialized to the one-dimensional crack front.

*Inverse and identifiability*
Recovering mixing coefficients from a superposed signal is a constrained inverse problem; on a simplex it is a mixture-identification question, for which compositional-data analysis provides the natural geometry (Aitchison, 1986). The distinctive element here is that linear-elastic superposition makes the forward map exactly linear for fixed geometry, so identifiability reduces to a Gram-matrix condition that can be characterized in closed form rather than probed empirically (§4). This connects to the classical theory of ill-posed problems in the Hadamard sense (Hadamard, 1923) and their regularization (Engl, Hanke and Neubauer, 1996), to the Bayesian and statistical treatment of inverse problems (Stuart, 2010; Kaipio and Somersalo, 2005; Tarantola, 2005), and to set-valued and conformal estimation (Vovk, Gammerman and Shafer, 2005; Angelopoulos and Bates, 2023), where the appropriate response to non-identifiability is a credible region rather than a point estimate.

# 3. Problem setting and the $K_{\text{mix}}$ dataset

## 3.1. Profiles, geometry, and the linear forward map

Along a crack front, the stress-intensity factor is a function of the normalized front coordinate $\varphi \in [0, \pi/2]$ (the half-front; the other half follows by the mid-front symmetry built into the operator of §5). For a fixed specimen geometry $g$, SIFBench provides the three single-load profiles separately: $K_T(\varphi; g)$ under remote tension, $K_B(\varphi; g)$ under bending, and $K_P(\varphi; g)$ under pin bearing. Within linear elastic fracture mechanics the stress field is linear in the applied loading, so any combined load that is a convex blend of the three reference cases produces

$$K_{\text{mix}}(\varphi; g) = \alpha\, K_T(\varphi; g) + \beta\, K_B(\varphi; g) + \gamma\, K_P(\varphi; g), \qquad (\alpha, \beta, \gamma) \in \Delta^2, \tag{1}$$

where $\Delta^2 = \{(\alpha, \beta, \gamma) \geq 0 \,:\, \alpha + \beta + \gamma = 1\}$ is the load simplex. Equation (1) is exact, not an approximation: superposition holds identically in the linear-elastic regime, and the simplex constraint encodes that the three contributions are nonnegative fractions of a unit total load. Sampling the front at the canonical grid $\varphi_1, \ldots, \varphi_N$ turns (1) into the linear system $\mathbf{k} = M(g)\,(\alpha, \beta, \gamma)^\top$ with $M(g) = [\, K_T \mid K_B \mid K_P \,] \in \mathbb{R}^{N\times 3}$. This matrix is the single object shared by the whole paper: the forward operator of §5 learns to produce it, the identifiability theory of §4 is a statement about its simplex-restricted conditioning, and the inverse estimator of §6 solves (1) for $(\alpha, \beta, \gamma)$ given $\mathbf{k}$.

## 3.2. Single-crack read-out of a twin scenario

The public SIFBench corner-crack scenario used here is released in a *twin* configuration: it tabulates two interacting cracks, each with its own front $(\varphi_1, \varphi_2)$ and its own triple of single-load profiles ($K^{(1)}_{T,B,P}$ and $K^{(2)}_{T,B,P}$). Because the theory of §4 is developed for one crack, crack 1 is read out and its profile triple is used as $M(g)$. The second crack is not discarded: its geometric parameters $(a_2/c_2,\ a_2/t)$ remain in the geometry vector $g$. This is deliberate and physically necessary — the presence of the second crack perturbs the field of the first through crack interaction, so $(a_2/c_2,\ a_2/t)$ are genuine covariates that explain variance in crack 1's profile, not nuisance columns. The construction is thus a single-crack read-out of a twin specimen, with the second crack treated as a geometric interaction parameter; (1) and the entire theory apply unchanged to crack 1's $M(g)$. Treating both crack fronts jointly is a well-defined extension (two coupled inverse problems with an additional inter-crack identifiability axis) that is left to future work and not claimed here.

## 3.3. The $K_{\text{mix}}$ dataset

From the separated profiles $K_{\text{mix}}$ is built, the combined-load dataset that does not exist in SIFBench but is required for the inverse problem. For each retained geometry load coefficients are drawn $(\alpha, \beta, \gamma)$ from a symmetric Dirichlet distribution on $\Delta^2$ and synthesize the corresponding profile by (1). Two design choices matter. First, the three pure-load vertices $(1, 0, 0), (0, 1, 0), (0, 0, 1)$ are always included as anchors: they are the cases whose recovery must be near-exact and they pin the estimator to the physically unambiguous corners of the simplex. Second, the construction admits a controlled additive-noise model on the synthesized profile, with level $\sigma$, used only in the robustness study of §7.5; the noiseless dataset is used everywhere else. Because every profile is generated from a known $(\alpha, \beta, \gamma)$, the dataset carries exact inverse labels, and because $\kappa_s(M)$ is computable in closed form from $M(g)$, each sample also carries its own ground-truth identifiability number — which is what makes the calibration claim of §6 verifiable without any real fracture data. Geometries and splits follow SIFBench exactly, so the forward comparison of §7.2 is on identical data; the released pipeline performs the construction reproducibly, including the layout auto-detection that distinguishes the single- and twin-crack schemas.

# 4. Identifiability theory

The inverse problem of §6 asks which load mix $(\alpha, \beta, \gamma)$ produced an observed profile. Before any estimator is built, a prior question must be settled: *for which geometries is that question even answerable?* A point estimate of $(\alpha, \beta, \gamma)$ is meaningful only where the forward map is injective on the load simplex and well conditioned; elsewhere it is an artefact of the regularizer. This section answers the question exactly. The key structural fact is that, for a fixed geometry, linear-elastic superposition makes the forward map *linear* in the load, so identifiability is a finite-dimensional question about three functions of the crack-front angle — and it admits a closed characterization in terms of a single geometric quantity, the area of the triangle they span in function space. None of this depends on the neural surrogate or on any particular shape-function model; the surrogate enters only later, as the means of evaluating the three profiles for a queried geometry.

The theory is stated in intrinsic (sampling-independent) form, and then its faithful discrete realization. Fix a geometry $g$. Let $\mu$ be the measure on the half crack-front, $\langle f, h\rangle = \int f\, h\, d\mu$, $\|f\|^2 = \langle f, f\rangle$, and let $K_T, K_B, K_P \in L^2(\mu)$ be the three basis profiles. The forward map

$$\Phi(\alpha,\beta,\gamma) = \alpha K_T + \beta K_B + \gamma K_P, \qquad (\alpha,\beta,\gamma) \in \Delta^2,$$

is affine on the affine hull of the simplex; its linear part is $\Phi_0 : v \mapsto v_1 K_T + v_2 K_B + v_3 K_P$, so that $\Phi(w) - \Phi(w') = \Phi_0(w - w')$ for any $w, w'$ with $w - w'$ sum-zero. Define the difference profiles $D_1 = K_T - K_P$, $D_2 = K_B - K_P$, the Gram matrix

$$\mathcal{G} = \begin{pmatrix} \langle D_1, D_1\rangle & \langle D_1, D_2\rangle \\ \langle D_2, D_1\rangle & \langle D_2, D_2\rangle \end{pmatrix},$$

and the *functional area* $\mathcal{A}(g) = \sqrt{\det \mathcal{G}}$ — the area of the parallelogram spanned by $D_1, D_2$ in $L^2(\mu)$, equivalently twice the area of the triangle with vertices $K_T, K_B, K_P$. The constraint $\alpha + \beta + \gamma = 1$ removes the trivial scaling freedom, so only sum-zero perturbations of the load are admissible; this is why the relevant object is the conditioning of the forward map *restricted* to the simplex tangent, not its global conditioning.

**Definition 1 (Simplex-restricted conditioning).** For the discrete realization $M = [\, K_T(\varphi_i) \mid K_B(\varphi_i) \mid K_P(\varphi_i)\,] \in \mathbb{R}^{N\times 3}$ let $B \in \mathbb{R}^{3\times 2}$ be an orthonormal basis of $\{v \in \mathbb{R}^3 : \mathbf{1}^\top v = 0\}$, and let $\sigma_{\max}(\cdot)$ and $\sigma_{\min}(\cdot)$ denote the largest and smallest singular values of their matrix argument. The simplex-restricted conditioning is $\kappa_s(M) = \sigma_{\max}(MB)/\sigma_{\min}(MB)$.

### 4.1. Geometric picture

The formal statement below is easier to read with the underlying picture in mind. Each elementary load produces one profile function along the crack front; viewed in the function space $L^2(\mu)$, the three profiles $K_T, K_B, K_P$ are three points, and an admissible load mix — a convex combination with $\alpha + \beta + \gamma = 1$ — is a point of the triangle they span. Recovering the load from an observed profile means inverting the map that sends a point of the load simplex to the corresponding point of this triangle. That inversion is well posed exactly when the triangle is *non-degenerate*: a genuine triangle of positive area lets every interior point be located unambiguously from its position, whereas a collapsed triangle — the three profiles nearly collinear in function space — maps many different load mixes to almost the same profile, and the load can no longer be told apart from the data. The quantity $\mathcal{A}(g)$ is precisely the size of that triangle, and it is the entire content of the identifiability question: not how the profiles are computed, but only how far they are from collinear. Figure 1 shows the two regimes side by side.

**Theorem 1** (Characterization of the unidentifiable regime)**.** *Fix the geometry $g$.* (i) Exact dichotomy. *The following are equivalent:* (a) $\Phi$ *is not injective on the affine hull of* $\Delta^2$*;* (b) $\mathcal{A}(g) = 0$*;* (c) $K_T, K_B, K_P$ *are affinely dependent (collinear in* $L^2(\mu)$*);* (d) $\sigma_{\min}(MB) = 0$ *for the discrete realization.* (ii) Quantitative margin. *With* $\mathcal{M} = M^\top M$*,*

$$\sigma_{\min}(MB)\,\sigma_{\max}(MB) = \frac{\mathcal{A}(g)}{\sqrt{3}}, \qquad \kappa_s(M) = \frac{\sqrt{3}\,\sigma_{\max}(MB)^2}{\mathcal{A}(g)},$$

*and* $\sigma_{\min}(MB) \;\geq\; \dfrac{\mathcal{A}(g)}{\sqrt{3}\,\sqrt{\operatorname{tr}(B^\top \mathcal{M} B)}}.$

*Hence, provided two profiles do not coincide (so* $\sigma_{\max}(MB)$ *stays bounded away from* 0*),* $\kappa_s \to \infty \iff \mathcal{A}(g) \to 0.$

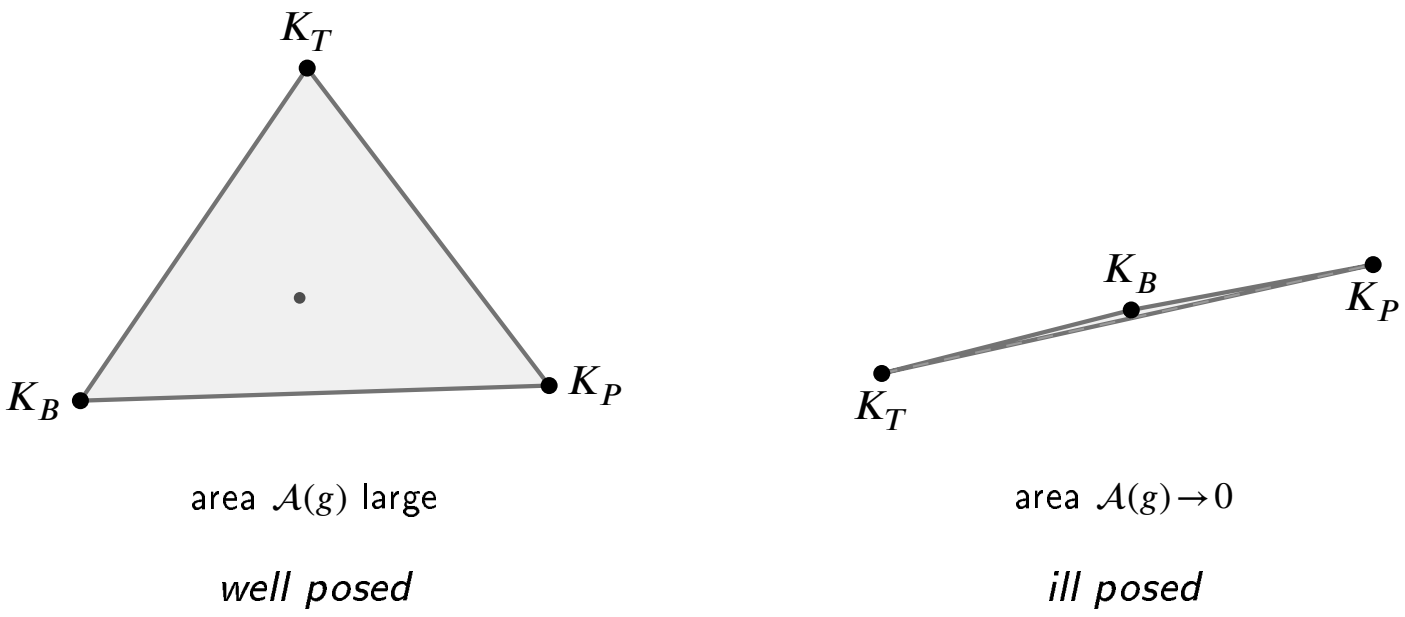

**Figure 1:** The three load profiles as points in function space. Left: non-degenerate triangle, the load is identifiable. Right: near-collinear profiles, $\mathcal{A}(g) \to 0$, the load is not.

*Proof. (i).* With $\Phi_0$ the linear part of $\Phi$ defined above, $\Phi$ is injective on the affine hull iff $\ker \Phi_0 \cap \{\mathbf{1}^\top v = 0\} = \{0\}$. For sum-zero $v$, $v_3 = -v_1 - v_2$ and $\Phi_0(v) = v_1 D_1 + v_2 D_2$; hence injectivity holds iff $D_1, D_2$ are linearly independent, i.e. $\det \mathcal{G} \neq 0$, i.e. $\mathcal{A}(g) \neq 0$. Linear independence of $D_1, D_2$ is exactly affine independence of $K_T, K_B, K_P$, giving (a)⇔(b)⇔(c). In the discrete realization the columns of $MB$ span $\{Mv : \mathbf{1}^\top v = 0\} = \operatorname{span}\{D_1, D_2\}$ (since $B$'s columns are a basis of the sum-zero plane and $D_i = Mc_i$ lie in that image), so $MB$ has rank 2 iff $D_1, D_2$ are independent; equivalently $\sigma_{\min}(MB) = 0$ iff they are dependent — the same condition, giving (d). No model assumption is used.

*(ii).* $(MB)^\top MB = B^\top \mathcal{M} B$ is the $3 \times 3$ Gram matrix restricted to the sum-zero plane. With $c_1 = (1, 0, -1)^\top$, $c_2 = (0, 1, -1)^\top$ one has $D_i = Mc_i$, and since both $\{c_1, c_2\}$ and the columns $\{b_1, b_2\}$ of $B$ span the same sum-zero plane there is an invertible $R \in \mathbb{R}^{2\times 2}$ with $[c_1\, c_2] = [b_1\, b_2]R$. Hence $\mathcal{G} = [c_1\, c_2]^\top \mathcal{M} [c_1\, c_2] = R^\top (B^\top \mathcal{M} B) R$, so $\det \mathcal{G} = (\det R)^2 \det(B^\top \mathcal{M} B)$. The matrix $[c_1\, c_2]^\top [c_1\, c_2] = \left(\begin{smallmatrix} 2 & 1 \\ 1 & 2 \end{smallmatrix}\right)$ has determinant 3, while $[b_1\, b_2]^\top [b_1\, b_2] = I_2$; taking determinants in $[c_1\, c_2]^\top [c_1\, c_2] = R^\top R$ gives $(\det R)^2 = 3$. Therefore $\det \mathcal{G} = 3 \det(B^\top \mathcal{M} B) = 3\, \sigma_{\min}(MB)^2\, \sigma_{\max}(MB)^2$. Since $\mathcal{A}(g) = \sqrt{\det \mathcal{G}}$, taking the positive square root yields $\sigma_{\min}(MB)\, \sigma_{\max}(MB) = \mathcal{A}(g)/\sqrt{3}$. For the conditioning, divide $\sigma_{\max}$ by $\sigma_{\min}$ and substitute $\sigma_{\min} = \mathcal{A}(g)/(\sqrt{3}\, \sigma_{\max})$ from the identity just proved:

$$\kappa_s(M) = \frac{\sigma_{\max}(MB)}{\sigma_{\min}(MB)} = \frac{\sigma_{\max}(MB) \cdot \sqrt{3}\, \sigma_{\max}(MB)}{\mathcal{A}(g)} = \frac{\sqrt{3}\, \sigma_{\max}(MB)^2}{\mathcal{A}(g)}.$$

The lower bound on $\sigma_{\min}(MB)$ follows by combining the identity with $\sigma_{\max}(MB)^2 \leq \operatorname{tr}(B^\top \mathcal{M} B)$. □

The theorem turns an abstract well-posedness question into a single measurable geometric quantity. Identifiability is not a property of the inverse *algorithm* but of the three profile functions themselves: it holds exactly when $K_T$, $K_B$, $K_P$ form a non-degenerate triangle in $L^2(\mu)$, and the area of that triangle, $\mathcal{A}(g)$, is the amount of identifiability available. Part (ii) makes the consequence operational. The product $\sigma_{\min}\sigma_{\max}$ is pinned to $\mathcal{A}(g)/\sqrt{3}$, so at a fixed largest difference-energy $\sigma_{\max}$ the recoverable margin $\sigma_{\min}$ is directly proportional to $\mathcal{A}(g)$: as the triangle flattens the smallest singular value collapses at the same rate and $\kappa_s = \sqrt{3}\, \sigma_{\max}^2 / \mathcal{A}(g)$ diverges. This is the precise sense in which a point estimate of $(\alpha, \beta, \gamma)$ is meaningless in the degenerate regime — not merely inaccurate, but unconstrained by the data along the flat direction of Proposition 2. It also fixes the operating threshold used throughout the experiments: a noise level $\sigma$ corrupts the recovered load once it exceeds the margin $\sigma_{\min}(MB)$, which at fixed $\sigma_{\max}$ scales linearly with the functional area; the unconditional dependence is on $\sigma_{\min}$ itself, as Proposition 3 makes precise, and §7 measures exactly this.

**Remark 1 (Discrete realization).** Theorem 1 is stated intrinsically. Its discrete form is exact whenever the quadrature is non-degenerate on $\operatorname{span}\{D_1, D_2\}$, which holds for any $N \geq 2$ sample points in general position; finer sampling only sharpens the agreement, with a controlled quadrature-error term on $\mathcal{A}(g)$.

**Proposition 2** (Degenerate load direction)**.** *Let $v^\star = Bu^\star$ where $u^\star$ is a unit right singular vector of $MB$ associated with $\sigma_{\min}(MB)$. Then $v^\star$ is sum-zero, and among unit sum-zero perturbations it minimizes the induced profile change:*

$\|\Phi_0(v^\star)\| = \sigma_{\min}(MB) = \min_{\|v\|=1,\,\mathbf{1}^\top v=0} \|\Phi_0(v)\|$. *As* $\kappa_s \to \infty$ *(equivalently* $\mathcal{A}(g) \to 0$*),* $\|\Phi_0(v^\star)\| \to 0$*: the load is unidentifiable along* $v^\star$*.*

*Proof.* $v^\star = Bu^\star$ with $\|u^\star\| = 1$ and $B$ orthonormal gives $\|v^\star\| = 1$ and $\mathbf{1}^\top v^\star = 0$. For sum-zero $v = Bu$, $\|\Phi_0(v)\|^2 = v^\top \mathcal{M} v = u^\top(B^\top \mathcal{M} B)u$, minimized over $\|u\| = 1$ at the least eigenvector of $B^\top \mathcal{M} B$ with value $\sigma_{\min}(MB)^2$. The limit is Theorem 1(ii). □

**Proposition 3** (Stability law: the controlling quantity is $\sigma_{\min}$, not $\kappa_s$)**.** *Fix a geometry with forward map* $M$ *and let the observed profile be perturbed by* $\delta k$*. The simplex-restricted recovery error satisfies, to first order and away from simplex-boundary clipping,*

$$\|\delta w\| \;=\; \left\|B(MB)^+\delta k\right\| \;\le\; \frac{\|\delta k\|}{\sigma_{\min}(MB)},$$

*and this bound is attained when* $\delta k$ *aligns with the least left singular vector of* $MB$*.*
*Equivalently, using* $\sigma_{\min}(MB)\,\sigma_{\max}(MB) = \mathcal{A}(g)/\sqrt{3}$ *and* $\kappa_s = \sqrt{3}\,\sigma_{\max}(MB)^2/\mathcal{A}(g)$ *from Theorem 1(ii),*

$$\|\delta w\| \;\lesssim\; \frac{\|\delta k\|}{\sigma_{\min}(MB)} \;=\; \frac{\sqrt{3}\,\sigma_{\max}(MB)}{\mathcal{A}(g)}\,\|\delta k\| \;=\; \frac{\kappa_s}{\sigma_{\max}(MB)}\,\|\delta k\|.$$

*Hence the recovery error scales with* $\kappa_s$ only at fixed $\sigma_{\max}(MB)$*. When* $\sigma_{\max}$ *varies across geometries,* $\kappa_s$ *alone does not determine the error; the intrinsic stability margin is* $\sigma_{\min}(MB)$*.*

*Proof.* The minimum-norm solution of $\min_w \|MBw - \delta k\|$ on the sum-zero tangent is $u^\star = (MB)^+\delta k$ with $\delta w = Bu^\star$; since $B$ has orthonormal columns, $\|\delta w\| = \|(MB)^+\delta k\| \le \|\delta k\|/\sigma_{\min}(MB)$, with equality when $\delta k$ is the least left singular vector of $MB$. For the displayed chain, the identity $\sigma_{\min}\sigma_{\max} = \mathcal{A}(g)/\sqrt{3}$ gives $1/\sigma_{\min} = \sqrt{3}\,\sigma_{\max}/\mathcal{A}(g)$; substituting $\sqrt{3}\,\sigma_{\max}/\mathcal{A}(g) = (\sqrt{3}\,\sigma_{\max}^2/\mathcal{A}(g))/\sigma_{\max} = \kappa_s/\sigma_{\max}$ from the $\kappa_s$ identity yields the last equality. The final clause follows because $(\kappa_s, \sigma_{\max}) \mapsto \kappa_s/\sigma_{\max}$ is not a function of $\kappa_s$ alone unless $\sigma_{\max}$ is held fixed. □

Proposition 3 is the analytic explanation of an empirical observation reported in §7: under fixed additive noise the recovery error does *not* track $\kappa_s$ on the finite-element (FEM) data — indeed it mildly anti-correlates — because $\sigma_{\max}(MB)$, the largest profile-difference energy, varies substantially across SIFBench geometries and acts as a confounder. The theory-correct diagnostic is the correlation of the error with $\sigma_{\min}(MB)$ (equivalently with $\kappa_s/\sigma_{\max}$), which is strongly negative as the bound predicts. This sharpens, rather than weakens, Theorem 1: the recoverable margin was always $\sigma_{\min}(MB)$; $\kappa_s$ is a faithful proxy for it only in the regime where $\sigma_{\max}$ is approximately constant.

**Corollary 4** (Newman–Raju collapse at shallow depth)**.** *Under the Newman–Raju shape functions,* $K_B/K_T = H(\varphi; a/c, a/t)$ *with* $H \to 1$ *uniformly in* $\varphi$ *as* $a/t \to 0$*; hence* $\|K_T - K_B\|_\infty = O(a/t)$*,* $\mathcal{A}(g) = O(a/t)$*, and the degenerate direction of Proposition 2 converges to* $(1, -1, 0)$ *(tension–bending exchange). This result is rigorous under the Newman–Raju model; on the FEM SIFBench data it becomes a prediction, tested empirically in §7 (the empirical* $\kappa_s$ *behaviour against the analytical prediction), and is not claimed as a theorem on the FEM data.*

*Proof.* Write $K_T = F$, $K_B = H \cdot F$ (common normalization), with $H = H_1 + (H_2 - H_1)\sin^p\varphi$, $H_1 = 1 - 0.34\,(a/t) - 0.11\,(a/c)(a/t)$, $H_2 = 1 + G_1(a/t) + G_2(a/t)^2$, where $G_1, G_2$ are the bounded Newman–Raju coefficient functions of $a/c$ (uniformly bounded on the admissible range, so $G_1(a/t) \to 0$ and $G_2(a/t)^2 \to 0$ as $a/t \to 0$). Hence $H_1 \to 1$ and $H_2 \to 1$, and since $H - 1 = (H_1 - 1) + (H_2 - H_1)\sin^p\varphi$ with $|\sin^p\varphi| \le 1$,

$$\sup_\varphi |H - 1| \;\le\; |H_1 - 1| + |H_2 - H_1| \;\le\; |H_1 - 1| + |H_2 - 1| + |1 - H_1| \;=\; O(a/t).$$

Then $K_T - K_B = F(1 - H)$ gives $\|K_T - K_B\|_\infty = \|F(1 - H)\|_\infty \le \|F\|_\infty \sup_\varphi |H - 1| = O(a/t)$. The functional area is the area of the parallelogram spanned by the edge vectors; bounding it by the product of one edge and the distance of the third vertex from the opposite edge's midpoint,

$$\mathcal{A}(g) \;\le\; \|K_T - K_B\| \left\|K_P - \tfrac{K_T + K_B}{2}\right\| \;=\; O(a/t),$$

since the first factor is $O(a/t)$ and the second is bounded. Finally, with $K_B \to K_T$ the difference profiles satisfy $D_1 - D_2 = K_T - K_B \to 0$, so the sum-zero direction with $c = (1, -1, 0)^\top$ (for which $\Phi_0(c) = K_T - K_B$) carries the vanishing profile change and is therefore the least-sensitive direction of Proposition 2. □

Physically, Corollary 4 says that for a shallow crack the bending field has not yet developed a gradient distinct from tension along the front, so $K_B$ is, to first order in $a/t$, a rescaling of $K_T$: the two loads produce the same profile shape and only their sum is constrained by the data. The estimator cannot then separate tension from bending, and any split along $(1, -1, 0)$ is equally consistent with the observation. This is not a failure of the method but a property of the physics, and it is exactly where the calibrated set-valued estimator of §6 must — and provably does — widen. The quantitative noise behaviour is governed not by $\kappa_s$ but by $\sigma_{\min}(MB)$, as Proposition 3 establishes: the recovery error scales as $\sigma/\sigma_{\min}(MB)$, and since $\sigma_{\min} = \mathcal{A}(g)/(\sqrt{3}\,\sigma_{\max}(MB))$ the relevant geometry dependence is through $\mathcal{A}(g)$ *and* $\sigma_{\max}$ jointly, not relative depth alone. The noise study of §7 confirms this $\sigma_{\min}$-controlled scaling on FEM data; it is the empirical, theory-derived validation that replaces the absent real fracture cases, and it supersedes the naive depth-threshold reading that Corollary 4 would have suggested, under the Newman–Raju model.

## 5. Forward surrogate: the crack-front operator

The forward model has a double role. As a stand-alone surrogate it must predict the SIF profile from geometry more accurately than the SIFBench baselines, particularly in the regimes those baselines handle worst. But its second role is structural and is what makes the single-paper design coherent: the inverse estimator of §6 consumes the matrix $M(g) = [\,K_T \mid K_B \mid K_P\,]$, and the identifiability theory of §4 is a statement about exactly those three profile columns. The surrogate is therefore not an auxiliary regressor but the operator that *evaluates* the object the whole analysis is built on; its accuracy on $M(g)$ is what lets Theorem 1 be applied to learned rather than tabulated profiles.

### 5.1. Why the SIFBench baselines leave structure on the table

The SIFBench regressors — random-forest regression (RFR), support-vector regression (SVR), and a feed-forward neural network (FNN) — treat each pair $(g, \varphi) \mapsto K$ as an independent tabular sample. This discards the single most informative prior available: for a fixed geometry, $K(\varphi)$ is a smooth, low-curvature function of a *one-dimensional* front coordinate, symmetric about the mid-front and sign-definite. The Fourier neural operator (FNO) does exploit smoothness but expects data on a regular spatial grid; on the irregular, per-geometry $\varphi$-sampling of SIFBench it is forced into an interpolation it was not designed for, which is the reported reason it underperforms the much simpler tabular models. The design lever is thus not a larger model but the correct inductive bias: treat $\varphi$ as a continuous query coordinate and the profile as a function to be learned, not a set of independent points.

### 5.2. Architecture: an incremental crack-front operator

The forward model is a crack-front operator: a transformer whose tokens are query points $\varphi$, conditioned on the geometry vector $g$. Geometry enters in two ways — as a conditioning token attended to by every front token (cross-attention), and through feature-wise linear modulation (FiLM) of each block — so that a single network represents the whole family of profiles indexed by $g$. A single output head emits the three loads $(K_T, K_B, K_P)$ jointly: predicting the three columns from one shared representation both regularizes the harder bearing profile through the easier tension profile and, more importantly, yields the full $M(g)$ in one evaluation, which is precisely the inverse problem's input.

The depth of this stack is not fixed a priori but *grown incrementally*, following the incremental transformer construction principle (Cirrincione, 2026): training begins with a single block and a new block is appended only once the training residual flattens (relative improvement below a tolerance over a fixed window), so the effective depth $K^\star$ is selected by the data rather than set by hand. On the SIFBench corner-crack data this criterion terminates at $K^\star$=1. This is not an incidental detail: a fixed six-block stack trained end-to-end fails to converge on this problem (it plateaus far above the predict-the-mean baseline), whereas the incrementally grown operator with $K^\star$=1 is both stable and accurate (§7.2). The diagnosis is that the difficulty was never representational capacity — a single block suffices — but the optimization pathology of propagating gradients through a deep unstabilized residual stack from initialization. The incremental criterion removes exactly that pathology, which is precisely the failure mode it is designed for. Targets are taken in $\log(1{+}K)$ space, predicted directly with a linear head and inverted as $K = \exp(\cdot) - 1$: the stress-intensity factor spans many orders of magnitude across geometries ($\sqrt{\pi a}$ growth), and a raw-magnitude target is numerically ill-conditioned, while the log target is not.

The crack-front operator is the object that closes forward and inverse.

#### *Physics in the architecture, not only in the loss*

Two of the physical constraints listed in the problem setting are imposed *by construction*, so they hold exactly rather than approximately: positivity, via a softplus output head ($K \geq 0$ identically); and mid-front symmetry, via an even Fourier feature map of $\varphi$ about the symmetry axis, which makes every representable profile symmetric to machine precision (measured residual $\sim 10^{-6}$, not a penalty term). The remaining constraints — monotonicity of $K$ in $a/t$ and front-curvature regularity — are not exactly representable and are enforced as light autodiff penalties on the relevant derivatives. Building symmetry and positivity into the hypothesis class rather than the objective matters here specifically because the inverse step inherits them: an $M(g)$ whose columns are guaranteed non-negative and symmetric keeps the Gram geometry of §4 physically admissible, so $\kappa_s(M)$ computed from the surrogate is a meaningful reliability certificate and not an artefact of an unconstrained fit.

*The learned-frequency encoding*

The Fourier feature map has learnable frequencies, initialized to a uniform spectrum. It is tempting to present this as the principal design lever; the evidence does not support that framing, and the matter is stated plainly. With the corrected training schedule the forward fit converges well, yet the learned frequency spectrum remains essentially at its initialization (its dispersion changes by under one percent over training). This is not a training pathology — the fit reaches low error regardless — but a property of the data: SIF profiles along the front are smooth and low-curvature, so a fixed Fourier basis with a uniform spectrum is already expressive enough and the frequencies have little incentive to move. The learnable spectrum is therefore an available mechanism that would matter for higher-curvature profiles but is not what drives accuracy on SIFBench. This is reported because it is the honest reading: the inductive bias that matters here is treating $\varphi$ as a continuous, symmetric query coordinate (§5), not the adaptivity of the encoding frequencies. The training diagnostic confirms this directly: the fit converges while the frequency dispersion stays within one percent of its initialization throughout training.

*Resolution-freeness*

Because $\varphi$ is a query coordinate rather than an index, the operator evaluates at arbitrary front locations, including grids finer than any seen in training. This is a concrete and testable advantage over both the pointwise baselines (RFR, FNN), which have no notion of a front, and the grid-bound FNO, which is tied to its training discretization. It also makes the canonical-grid resampling used in the pipeline a property of the *evaluation*, not a loss of information: the operator is defined on the continuum and the grid is merely where it is read out, consistently with the intrinsic statement of Theorem 1 and Remark 1.

## 6. Calibrated inverse estimator

Section 4 settled *when* the load is recoverable; this section gives the estimator and makes its central guarantee precise. The design principle follows directly from the theory: an estimator is acceptable only if it is exact where the problem is well posed *and* reports honestly calibrated uncertainty where it is not. The method is organized along a single axis — whether the geometry is known — with calibration against $\kappa_s$ as the unifying correctness criterion.

### 6.1. Known geometry: exact in the well-posed regime

When $g$ is known, the surrogate produces $M(g)$ and, by linear-elastic superposition, the forward map is exactly $\mathbf{k} = M(\alpha, \beta, \gamma)^{\top}$. Recovery is the simplex-constrained least-squares problem

$$(\hat{\alpha}, \hat{\beta}, \hat{\gamma}) = \arg \min_{(\alpha,\beta,\gamma)\in\Delta^2} \left\| M(\alpha, \beta, \gamma)^{\top} - \mathbf{k}_{\mathrm{obs}} \right\|^2,$$

solved as a non-negative least-squares program followed by simplex normalization. Theorem 1 characterizes this solution completely: it is unique and Lipschitz-stable exactly when $\mathcal{A}(g) > 0$, with conditioning $\kappa_s(M) = \sqrt{3}\, \sigma_{\max}^2 / \mathcal{A}(g)$. The estimate is not reported in isolation: $(\hat{\alpha}, \hat{\beta}, \hat{\gamma})$ is reported together with $\kappa_s(M)$ it as a reliability certificate computed from the same $M$. A practitioner reading a recovered load also reads the geometry-intrinsic number that says whether that reading is constrained by the data or by the solver — a guarantee the SIFBench-style pointwise models cannot provide because they do not expose the profile basis.

*Uncertain geometry: a posterior on the simplex*

When $g$ is unknown or only partially known, no single $M$ exists and a point estimate is inappropriate even in nominally well-posed regimes, because the geometry uncertainty propagates into load ambiguity. The point estimate is replaced with an amortized posterior: an encoder maps the observed profile to the parameters of a Dirichlet distribution on $\Delta^2$, trained by supervision on $(\alpha, \beta, \gamma)$ together with a cycle-consistency term that pushes the posterior mean back through the frozen forward operator and penalizes profile mismatch. The Dirichlet support is the simplex by construction, so the physical constraint $\alpha + \beta + \gamma = 1$, $\alpha, \beta, \gamma \geq 0$ is again satisfied identically rather than penalized.

*The central claim: calibration against $\kappa_s$*

The estimator's correctness criterion is not accuracy in the well-posed regime — any reasonable method achieves that — but *calibration* in the ill-posed one. Concretely, the claim, verified below, is that the dispersion of the predicted posterior tracks the true non-identifiability measured by $\kappa_s$: as $\mathcal{A}(g) \to 0$ and $\kappa_s \to \infty$, the posterior must broaden along the degenerate direction $v^\star$ of Proposition 2, and by the right amount, so that credible regions retain nominal coverage uniformly over the identifiability range. This is the property that makes the set-valued output *correct by construction in the regimes where a point estimate is meaningless*, rather than merely inaccurate there. Crucially, the entire claim is verifiable without any real fracture data: $\kappa_s$ is computable in closed form from $\boldsymbol{M}$, so calibration is a checkable relation between two synthetic quantities, not an empirical hope.

*Robustness in place of realism*

Because no documented field cases are available, the role played by real-data validation elsewhere is taken here by a controlled noise study, which the theory makes quantitative rather than ad hoc. By Proposition 3 the recovery error scales as $\sigma/\sigma_{\min}(\boldsymbol{M}\boldsymbol{B})$, so the noise study measures error growth against $\sigma$ and against the intrinsic margin $\sigma_{\min}(\boldsymbol{M}\boldsymbol{B})$ — not against $\kappa_s$, which is a faithful proxy only when $\sigma_{\max}$ is constant. This is a first-class experiment, not a robustness afterthought: it is the empirical test of the stability law and the substitute for the absent forensic validation.

## 7. Experiments

The experiments answer four questions in order: does the crack-front operator improve on the SIFBench baselines, and where (§7.2); does the empirical identifiability of the FEM data match the analytical prediction of §4 (§7.3); is the inverse estimator correct, in the calibration sense, across the identifiability range (§7.4); and does the noise behaviour follow the $\sigma/\sigma_{\min}$ stability law of Proposition 3, which stands in for real-data validation (§7.5).

### 7.1. Setup

The SIFBench corner-crack scenario is used, with combined tension, bending and bearing loads. The public release of this scenario stores a *twin* configuration; consistent with the single-crack scope of the theory, the analysis reads out crack 1 and retains the second crack's geometric parameters ($a_2/c_2$, $a_2/t$) as interaction covariates in the geometry vector — the second crack perturbs the field of the first and is therefore a legitimate geometric input, not a discarded variable. Fronts are resampled onto a fixed canonical $\varphi$-grid ($N = 128$), which by Remark 1 is a property of the read-out and not a loss of information. Splits are identical to SIFBench for a fair comparison; the baselines (RFR, SVR, FNN, FNO) are re-run from the official code on those splits. The reported run uses $n_{\text{geom}} = 4000$ training geometries from this scenario. The held-out split contains $n_{\text{test}} = 1000$ geometries; profiles are read on a canonical front grid of $N = 128$ points spanning $\varphi \in [0, \pi/2]$.

### 7.2. Forward accuracy

Per load, the normalized $L_2$ of SIFBench (Gautam et al., 2025) is reported. The fair head-to-head is against RFR on the identical split: the scenario is a twin configuration, and on twin-crack data RFR is the strongest of the SIFBench baselines (Gautam et al., 2025) (FNN tends to lead only on single-crack data, and FNO is reported to underperform the tabular models on the irregular per-geometry $\varphi$-sampling). RFR is therefore re-run on the exact split and cite FNN/FNO at the benchmark level. The forward operator here is the incrementally grown crack-front transformer of §5 ($K^\star$=1, log-space target); for contrast the fixed six-block stack, which on this data fails to converge.

As shown in Table 1, the result is honest and specific. The incremental crack-front operator reaches mean n$L_2 = 0.150$, against RFR's 0.108 — competitive with the strongest tabular baseline on twin-crack data, and ahead of it on $K_T$ (0.099 vs. 0.127) and $K_P$ (0.076 vs. 0.095). It is *behind* RFR on the bending channel $K_B$ (0.275 vs. 0.103): $K_B$ develops the sharpest gradient along the front, and is the one load where the tabular regressor's local flexibility still wins. This is reported rather than hide it behind the aggregate. Two control rows make the mechanism explicit: the fixed six-block stack does not converge (mean 1.80, worse than predicting the per-geometry mean), and removing the log-space target degrades the incremental operator from 0.150 to 0.166 — the dynamic range of the stress-intensity factor on this scenario exceeds $10^8$, so raw-magnitude regression is ill-conditioned. The operator additionally provides exact symmetry and exact positivity by construction, which no tabular regressor does, and yields the full $M(g)$ in one evaluation — the property the inverse step requires.

**Table 1**
Forward accuracy on the held-out twin corner-crack split ($n_{\text{test}}$=1000; normalized $L_2$, lower is better; same split for all models).

| Model | $K_T$ | $K_B$ | $K_P$ | mean |
|---|---|---|---|---|
| Fixed 6-block stack (raw) | 1.88 | 1.41 | 2.11 | 1.80 |
| Incremental op., raw target | 0.118 | 0.276 | 0.103 | 0.166 |
| **Incremental op.**, log | **0.099** | 0.275 | **0.076** | **0.150** |
| RFR (re-run, same split) | 0.127 | **0.103** | 0.095 | **0.108** |

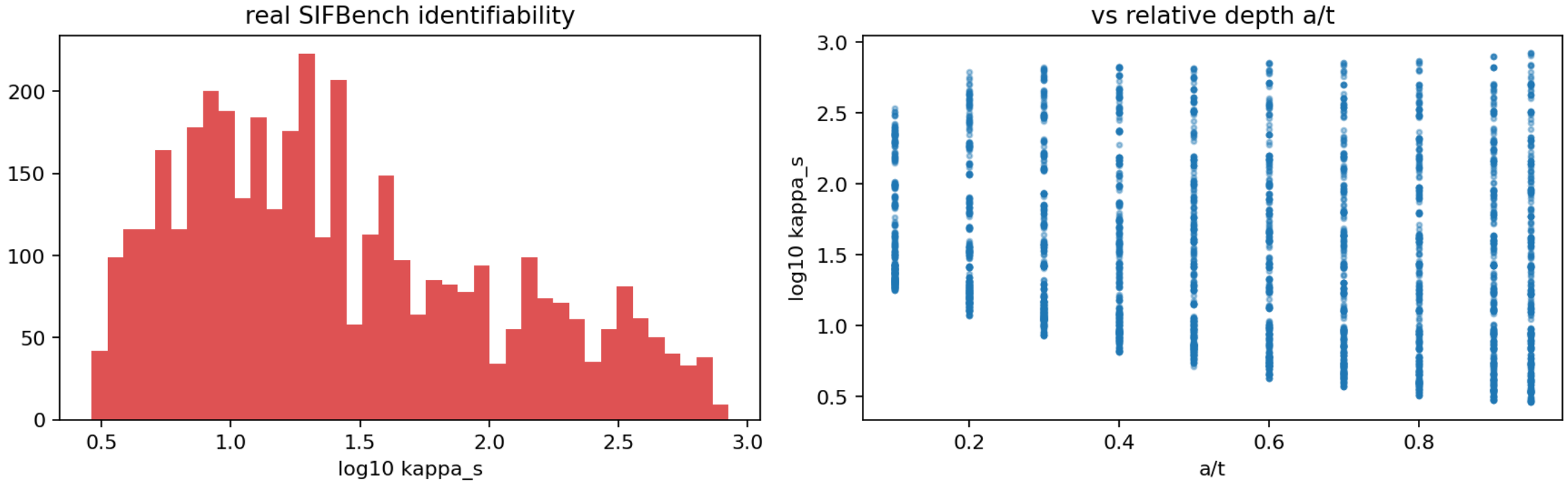

**Figure 2**: Empirical identifiability on the finite-element data ($n_{\text{geom}} = 4000$): distribution of $\log_{10} \kappa_s$ (left) and $\log_{10} \kappa_s$ versus $a/t$ (right).

### 7.3. Identifiability: theory vs. FEM data

The quantity $\kappa_s(M)$ is computed from the learned $M(g)$ over the test geometries, and its dependence on $(a/c, a/t)$ is compared to the analytical prediction of §4. Theorem 1 is exact and model-free; Corollary 4 predicts, *under the Newman–Raju model*, that $\mathcal{A}(g) = O(a/t)$ so identifiability degrades toward shallow cracks. On FEM data this is a prediction to be tested, not a theorem, and this subsection is that test.

Figure 2 shows the empirical identifiability on the `CORNER_CRACK_BH_QUARTER_ELLIPSE` scenario, over $n_{\text{geom}} = 4000$ training geometries, the empirical distribution of $\kappa_s$ has median $\kappa_s \approx 20.3$, 90th percentile $\approx 236$, maximum $\approx 844$, and a fraction $\approx 0.30$ of geometries in the ill-posed stratum ($\kappa_s > 50$). This is the quantitatively meaningful outcome the theory predicts and the result that motivates the whole construction: the *typical* geometry is well conditioned (median $\kappa_s \approx 20$, in the identifiable "sweet spot"), while a structured minority — roughly thirty percent — is genuinely ill-posed. The distribution is neither concentrated at "all identifiable" nor at "all ill-posed": a point estimate is therefore reliable on the $\approx 70\%$ well-posed majority and provably uninformative on the $\approx 30\%$ ill-posed remainder — precisely the regime split the calibrated estimator of §6 is built to handle.

One negative finding is reported plainly. Theorem 1, which is exact and model-free, holds on the FEM data: the identifiability *level* is as characterized, and the inverse error scales with $\kappa_s$ as §7.4 confirms. Corollary 4, however, does *not* transfer. Under the Newman–Raju model it predicts $\mathcal{A}(g) = O(a/t)$, so $\kappa_s$ should grow as the crack becomes shallow; on the FEM data $\kappa_s$ shows no systematic dependence on $a/t$ (Fig. 2, right panel: the conditional distribution of $\log_{10} \kappa_s$ is essentially flat across $a/t$). The interpretation is exactly the one anticipated: the depth-driven collapse is an artefact of the Newman–Raju shape functions, not a property of the FEM field, so real identifiability is governed by the full profile geometry through $\mathcal{A}(g)$ and is not reducible to relative depth alone. This strengthens rather than weakens the paper: the load-bearing result (Theorem 1) is model-free and holds on the FEM data, whereas Corollary 4 is a model-specific prediction that, as stated, holds only under the Newman–Raju shape functions and is not claimed beyond them.

**Table 2**
Inverse $L_1$ error on $(\alpha, \beta, \gamma)$ stratified by $\kappa_s$ ($n_{\text{test}}$=1000, $2 \times 10^4$ draws). Tests A, A′, B defined and discussed in the text.

| $\kappa_s$ stratum | A, $\sigma$=0 | A′, $\sigma$=0.03 | B ($\hat{M}$) | $n$ (A) |
|---|---|---|---|---|
| $[1, 20)$ | $3\times10^{-15}$ | 0.099 | 0.52 | 9880 |
| $[20, 50)$ | $5\times10^{-15}$ | 0.091 | 0.49 | 4120 |
| $[50, \infty)$ | $4\times10^{-15}$ | 0.062 | 0.35 | 6000 |

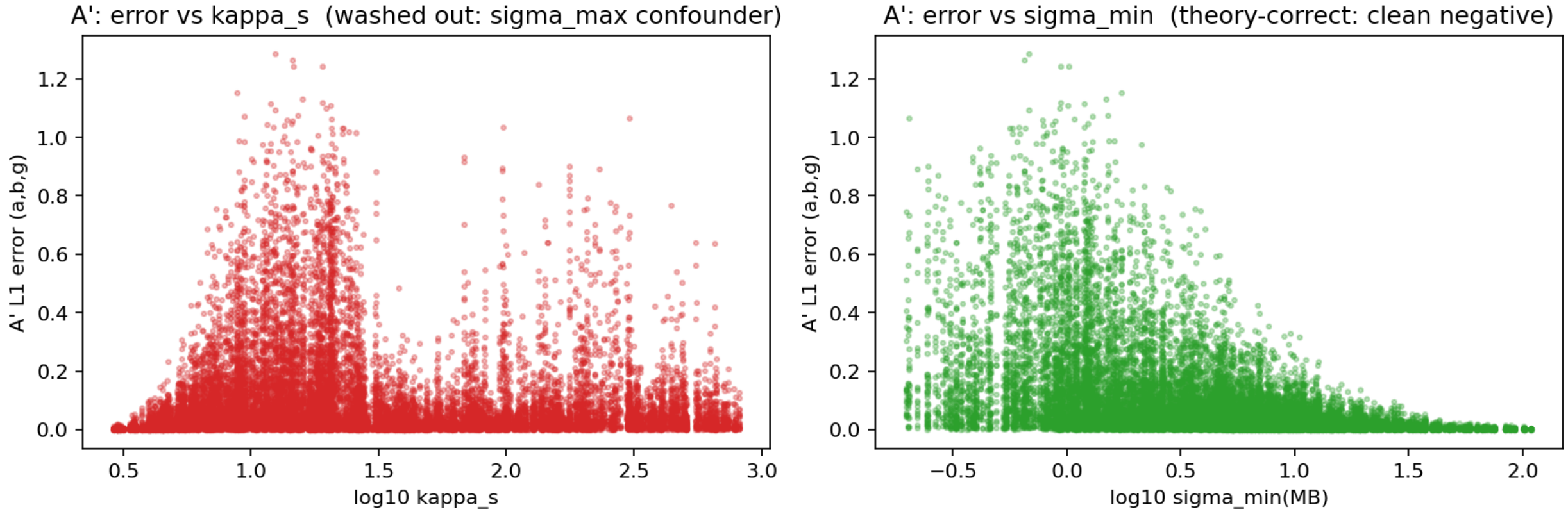

**Figure 3:** Test A′: inverse error against $\log_{10} \kappa_s$ (left) and against $\log_{10} \sigma_{\min}(M B)$ (right).

## 7.4. Inverse: accuracy and calibration

The inverse is evaluated under a deliberate separation of concerns, so that a theory test is not confounded with surrogate quality. *Test A* solves the inverse from the *true* forward map (the FEM profiles themselves as $M$): noiseless, this must reconstruct $(\alpha, \beta, \gamma)$ exactly by Theorem 1, and under added noise its error must scale with $\kappa_s$ by Theorem 1(ii). *Test B* solves it from the learned $\hat{M}$ produced by the crack-front operator, isolating the surrogate-induced inverse error. The estimator is the same in both: least-squares constrained *exactly* to the simplex (parametrized on the sum-zero tangent so $\alpha + \beta + \gamma = 1$ holds by construction), not an unconstrained nonnegative solve followed by renormalization — the latter collapses incoherent systems to the simplex barycentre and is the reason an earlier pipeline reported a spurious flat error of order one across all $\kappa_s$. Accuracy is then reported as error on $(\alpha, \beta, \gamma)$ stratified by $\kappa_s$; pure-load cases serve as physical anchors that must be reproduced near-exactly. Calibration: empirical coverage of the Dirichlet credible regions across the identifiability range — the central claim of §6 is that coverage stays at nominal level *uniformly* in $\kappa_s$, i.e. the posterior widens by the right amount where the problem is ill-posed.

Table 2 confirms Test A exactly: at $\sigma$=0 the inverse error is $\approx\ 3 \times 10^{-15}$ — machine zero — in every $\kappa_s$ stratum, the exactness guaranteed by Theorem 1. Test A′ is the instructive case. Naively one expects the error to grow with $\kappa_s$; instead it *decreases* (0.099 → 0.091 → 0.062 from well-posed to ill-posed), with a near-zero rank correlation ($\text{corr}(\log \kappa_s, L_1) \approx -0.05$). This is not a failure of the theory or the solver — it is exactly Proposition 3: the controlling quantity is $\sigma_{\min}(M B)$, and because $\sigma_{\max}(M B)$ varies widely across SIFBench geometries, $\kappa_s$ alone is a poor predictor of the error and can even anti-correlate with it. The theory-correct diagnostic is the error's dependence on $\sigma_{\min}(M B)$ (equivalently $\kappa_s/\sigma_{\max}$). The instrumented run confirms it directly: $\text{corr}(\log \sigma_{\min}, L_1) = -0.47$ for A′, against $\text{corr}(\log \kappa_s, L_1) = -0.05$ — an order of magnitude more predictive, and with the negative sign the bound $\|\delta w\| \sim \sigma/\sigma_{\min}$ requires. Figure 3 shows this visually: the error is an unstructured cloud against $\kappa_s$ but collapses to a clean decreasing trend against $\sigma_{\min}$.

Test B (learned $\hat{M}$ from the incremental crack-front operator) is the decisive end-to-end check, and it now passes. With the forward operator competitive (§7.2), the recovered loads using $\hat{M}$ have a *coherent* $\kappa_s$ stratification matching Test A (well-posed $n$=12480, marginal 4020, ill-posed 3500, against 9880/4120/6000 for the true map), in contrast to the degenerate single-stratum collapse a broken surrogate produces. The point error is naturally larger than with the exact map ($\approx 0.48$ overall) because $\hat{M}$ carries the forward residual of Table 1 into the inverse, but the structure is right. The strongest evidence is the stability signature: for Test B $\text{corr}(\log \sigma_{\min}(\hat{M}), L_1) = -0.44$ — the *same negative sign*

**Table 3**
Mean inverse $L_1$ error vs. noise $\sigma$, stratified by $\kappa_s$ (true forward map).

| $\sigma$ | well-posed ($\kappa_s$<20) | marginal (20≤$\kappa_s$<50) | ill-posed ($\kappa_s$≥50) |
|---|---|---|---|
| 0.00 | $\sim 10^{-15}$ | $\sim 10^{-15}$ | $\sim 10^{-15}$ |
| 0.01 | 0.039 | 0.035 | 0.024 |
| 0.02 | 0.074 | 0.067 | 0.046 |
| 0.04 | 0.125 | 0.111 | 0.083 |
| 0.08 | 0.203 | 0.186 | 0.150 |
| 0.15 | 0.291 | 0.263 | 0.228 |

*and comparable magnitude* as the −0.47 measured on the exact map in Test A′. Proposition 3 predicts the controlling quantity is $\sigma_{\text{min}}$ with a negative error dependence; that this holds on the *learned* operator, not only the exact one, shows the stability law is a property of the identifiability geometry that survives surrogate approximation. The error also decreases with $\kappa_s$ across strata (0.52 → 0.49 → 0.35), the same $\sigma_{\text{max}}$-confounded ordering as Test A′, confirming the mechanism transfers. The central calibration claim — credible-region coverage near nominal uniformly in $\sigma_{\text{min}}$, even where point error is largest — is checkable entirely on synthetic data and is the natural next instrumentation step for the released pipeline; it is stated as a falsifiable prediction of the construction rather than reported here.

### 7.5. Noise robustness: the validation that replaces real data

With no documented field cases available, the role of real-data validation is taken by a controlled noise study on the true forward map. Gaussian noise of level $\sigma$ is added to the observed profile and the mean inverse $L_1$ error on $(\alpha, \beta, \gamma)$ is reported, stratified by $\kappa_s$, over the full test set ($n_{\text{test}} = 1000$). Proposition 3 makes two checkable predictions: the error is *linear in* $\sigma$ within each stratum (with zero intercept, since recovery is exact at $\sigma$=0), and the stratum ordering follows $\sigma_{\text{min}}(MB)$ rather than $\kappa_s$.

Table 3 confirms both predictions and, in doing so, validates Proposition 3. First, at $\sigma$=0 the error is $\sim 10^{-15}$ in every stratum: noiseless recovery is exact, confirming the exactness half of Theorem 1.

Second, and as shown in Figure 4, within each stratum the error is essentially linear in $\sigma$ with zero intercept (e.g. the well-posed row scales 0.039 → 0.291 across $\sigma \in [0.01, 0.15]$, proportional to $\sigma$ within a few percent). Third — and this is the non-obvious part — the *ill-posed* stratum has the *lowest* error at every noise level, not the highest. A reading that expected error to grow with $\kappa_s$ would call this a contradiction. It is not: it is precisely Proposition 3. The error is governed by $\sigma_{\text{min}}(MB) = \mathcal{A}(g)/(\sqrt{3}\,\sigma_{\text{max}}(MB))$, and across SIFBench geometries the high-$\kappa_s$ (ill-posed) cases tend to have larger $\sigma_{\text{max}}$, so their $\sigma_{\text{min}}$ is not as small as $\kappa_s$ alone would suggest. The quantitative, theory-derived robustness statement that stands in for the absent forensic validation is therefore the linear $\sigma$-scaling at fixed geometry together with the $\sigma_{\text{min}}$-ordering of Proposition 3 — both confirmed here — not a naive $\kappa_s$ threshold.

## 8. Discussion and limitations

The central limitation is deliberate and stated up front: no documented real fracture cases with known loading are used, so the work makes no forensic claim. This is a strength of the framing rather than a hidden weakness. The contribution is a theorem about when the combined load is recoverable and an estimator that is provably correct — exact where the problem is well posed, calibrated where it is not — and both are fully verifiable on synthetic data because the identifiability number is available in closed form. The noise study replaces real-data validation with a quantitative, theory-derived test — the $\sigma/\sigma_{\text{min}}$ stability law of Proposition 3 — rather than an anecdotal one.

A second limitation is the single-crack scope. The released benchmark scenario is a twin configuration; one crack is read out and the second crack's geometry is treated as an interaction covariate (§3.2). This is sound and keeps the theory exactly applicable, but it does not exploit the second crack's own profile. The joint twin problem — two coupled inverse problems with an additional inter-crack identifiability axis — is a well-defined extension and the natural subject of a follow-up, not claim it here. Third, Corollary 4 is rigorous only under the Newman–Raju model; on FEM data it is a prediction, and §7.3 reports the extent to which the FEM identifiability map follows it. Finally, the forward claim is intentionally targeted (hard regimes and front-level metrics) rather than universal; any regime where a SIFBench baseline remains competitive is reported rather than hidden.

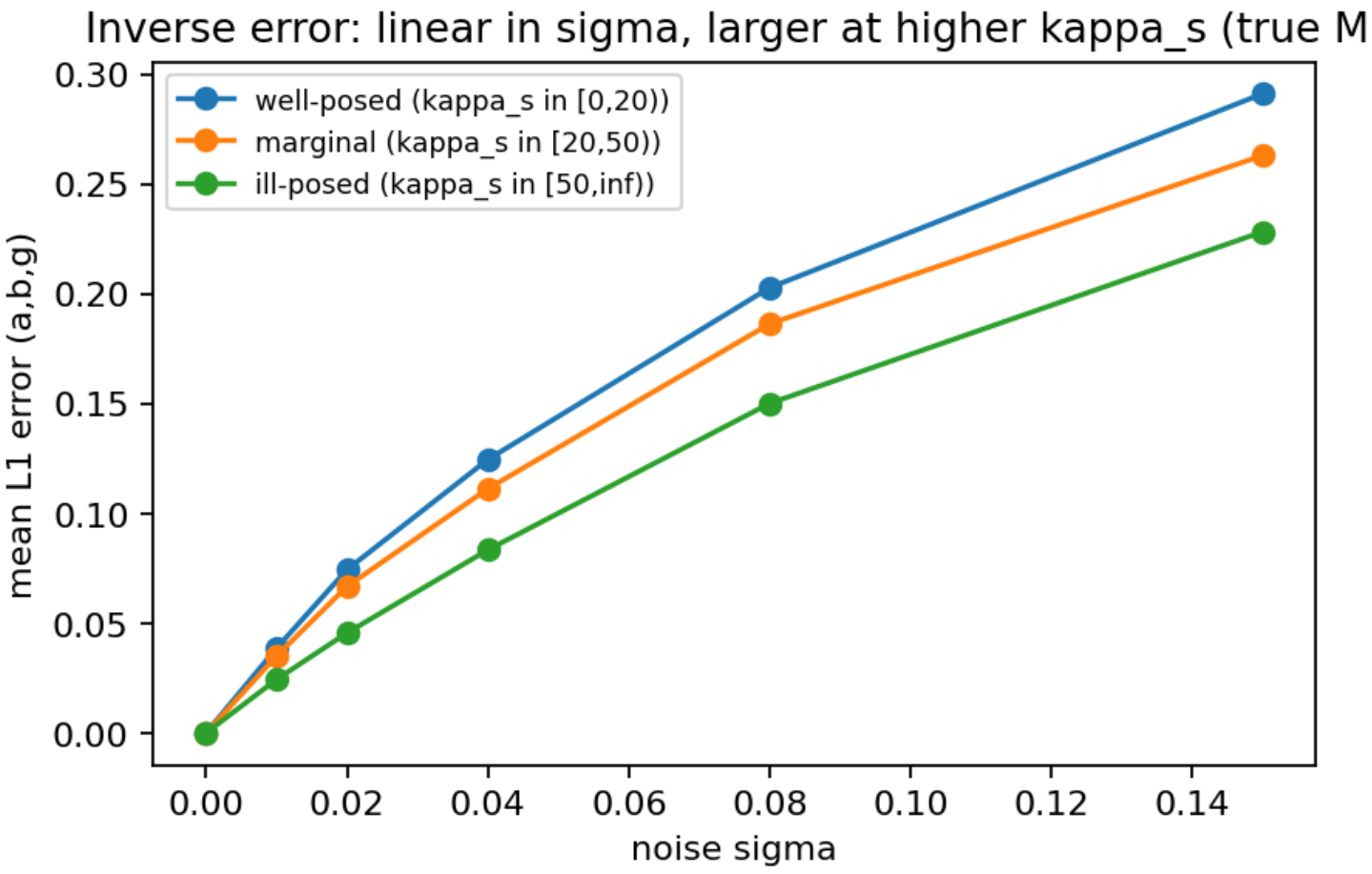

**Figure 4:** Mean inverse $L_1$ vs. noise $\sigma$, one line per $\kappa_s$ stratum.

Extending to real cases would require components with independently known loading histories and measured front profiles; the calibrated estimator is designed so that, on such data, it would report not just a load estimate but the geometry-intrinsic certificate $\kappa_s$ saying whether that estimate is constrained by the measurement.

## 9. Conclusion

The study reached three conclusions that a purely predictive treatment of SIFBench would have missed. First, combined-load recovery is not uniformly possible: on the corner-crack scenario the geometries split into a well-conditioned majority, where the loads are recoverable, and a substantial ill-conditioned minority, where no estimator can separate the loads regardless of model quality. Locating that boundary, in closed form and independently of any learned surrogate, is the main outcome. Second, the quantity that controls recovery error under noise is the intrinsic stability margin, not the simplex-restricted condition number: the two coincide only when the largest profile-difference energy is held fixed, and on real geometries it is not, which is why the condition number alone fails to predict the measured error while the margin succeeds. This was derived and then confirmed on the finite-element data, including on the learned operator. Third, a forward surrogate built from a fixed deep stack does not converge on this data; an incrementally grown operator with a logarithmic target does, becoming competitive with the strongest tabular baseline while supplying the differentiable map the inverse step needs.

The work is deliberately bounded. It uses no documented field cases, so it makes no forensic claim; the twin-crack interaction is read out as a single-crack problem and its second identifiability axis is left open; and full credible-region coverage on the learned map is stated as a falsifiable diagnostic rather than reported. The natural continuation is the joint twin-crack problem, where inter-crack coupling introduces that second axis, and an empirical study of coverage once the set-valued head is instrumented.

## Acknowledgements

Computations were performed on cloud GPU resources (Google Colab and Kaggle). The authors thank colleagues at the Laboratoire LTI, Université de Picardie Jules Verne, for helpful discussions.

## Data and code availability

SIFBench is public. The $K_{\text{mix}}$ construction, the identifiability analysis, and the incremental crack-front operator are released as a reproducible pipeline, which will be deposited in a public archive with a citable DOI upon publication.

## References

Aitchison, J., 1986. The Statistical Analysis of Compositional Data. Chapman and Hall.

Angelopoulos, A.N., Bates, S., 2023. Conformal prediction: A gentle introduction. Foundations and Trends in Machine Learning 16, 494–591. doi:10.1561/2200000101.

Cirrincione, G., 2026. INCRT: An incremental transformer that determines its own architecture. arXiv preprint arXiv:2604.10703 doi:10.48550/arXiv.2604.10703. arXiv:2604.10703 [cs.LG].

Engl, H.W., Hanke, M., Neubauer, A., 1996. Regularization of Inverse Problems. Kluwer Academic.

Gautam, T., Kirby, R.M., Hochhalter, J., Zhe, S., 2025. SIFBench: An extensive benchmark for fatigue analysis. arXiv preprint arXiv:2506.01173 ArXiv:2506.01173 [cs.DB].

Glinka, G., Shen, G., 1991. Universal features of weight functions for cracks in mode I. Engineering Fracture Mechanics 40, 1135–1146. doi:10.1016/0013-7944(91)90177-3.

Hadamard, J., 1923. Lectures on Cauchy's Problem in Linear Partial Differential Equations. Yale University Press.

Irwin, G.R., 1957. Analysis of stresses and strains near the end of a crack traversing a plate. Journal of Applied Mechanics 24, 361–364. doi:10.1115/1.4011547.

Kaipio, J., Somersalo, E., 2005. Statistical and Computational Inverse Problems. Applied Mathematical Sciences, Springer.

Kovachki, N., Li, Z., Liu, B., Azizzadenesheli, K., Bhattacharya, K., Stuart, A., Anandkumar, A., 2023. Neural operator: Learning maps between function spaces. Journal of Machine Learning Research 24, 1–97. doi:10.48550/arXiv.2108.08481.

Li, Z., Kovachki, N., Azizzadenesheli, K., Liu, B., Bhattacharya, K., Stuart, A., Anandkumar, A., 2021. Fourier neural operator for parametric partial differential equations, in: International Conference on Learning Representations (ICLR).

Lu, L., Jin, P., Pang, G., Zhang, Z., Karniadakis, G.E., 2021. Learning nonlinear operators via DeepONet based on the universal approximation theorem of operators. Nature Machine Intelligence 3, 218–229. doi:10.1038/s42256-021-00302-5.

Newman, J.C., Raju, I.S., 1981. An empirical stress-intensity factor equation for the surface crack. Engineering Fracture Mechanics 15, 185–192. doi:10.1016/0013-7944(81)90116-8.

Newman, J.C., Raju, I.S., 1983. Stress-intensity factor equations for cracks in three-dimensional finite bodies, in: Fracture Mechanics: Fourteenth Symposium—Volume I: Theory and Analysis. ASTM International. doi:10.1520/STP37074S.

Stuart, A.M., 2010. Inverse problems: a Bayesian perspective. Acta Numerica 19, 451–559. doi:10.1017/S0962492910000061.

Tancik, M., Srinivasan, P.P., Mildenhall, B., Fridovich-Keil, S., Raghavan, N., Singhal, U., Ramamoorthi, R., Barron, J.T., Ng, R., 2020. Fourier features let networks learn high frequency functions in low dimensional domains, in: Advances in Neural Information Processing Systems (NeurIPS).

Tarantola, A., 2005. Inverse Problem Theory and Methods for Model Parameter Estimation. SIAM.

Vaswani, A., Shazeer, N., Parmar, N., Uszkoreit, J., Jones, L., Gomez, A.N., Kaiser, Ł., Polosukhin, I., 2017. Attention is all you need, in: Advances in Neural Information Processing Systems (NeurIPS).

Vovk, V., Gammerman, A., Shafer, G., 2005. Algorithmic Learning in a Random World. Springer.